\documentclass[preprint,11pt,3p,authoryear]{elsarticle}

\usepackage[utf8]{inputenc}
\usepackage{amsmath}
\usepackage{amsfonts}
\usepackage{amsthm}
\usepackage{xcolor}
\usepackage{booktabs}
\usepackage{graphicx}
\usepackage[authoryear]{natbib}
\usepackage{bm}
\usepackage[ruled,vlined,linesnumbered]{algorithm2e}
\usepackage[bookmarks=false,colorlinks]{hyperref}
\AtBeginDocument{%
	\hypersetup{
		linkcolor=blue,   
		citecolor=blue,
		urlcolor=blue}}

\newcommand{\V}{\mathcal{V}}

\newcommand{\C}{\mathcal{C}}
\newcommand{\T}{\mathcal{T}}
\newcommand{\K}{\mathcal{K}}
\newcommand{\Ni}{\mathcal{N}_i}
\newcommand{\NGi}{\mathcal{N}^{\text{G}}_i}
\newcommand{\NAi}{\mathcal{N}^{\text{A}}_i}
\newcommand{\Vk}{\V_k}
\newcommand{\Vkp}{\V_{k'}}
\newcommand{\Vnk}{\V_{-k}}

\newcommand{\fij}{f_{ij}}
\newcommand{\pij}{p_{ij}}
\newcommand{\pji}{p_{ji}}

\newcommand{\sit}{S_i(t)}

\newcommand{\hsit}{\bar{S}_i(t)}
\newcommand{\hiit}{\bar{I}_i(t)}
\newcommand{\hrit}{\bar{R}_i(t)}
\newcommand{\hdit}{\bar{D}_i(t)}
\newcommand{\hsitp}{\bar{S}_i(t+1)}
\newcommand{\hiitp}{\bar{I}_i(t+1)}
\newcommand{\hritp}{\bar{R}_i(t+1)}
\newcommand{\hditp}{\bar{D}_i(t+1)}
\newcommand{\hsjt}{\bar{S}_j(t)}
\newcommand{\hijt}{\bar{I}_j(t)}
\newcommand{\hrjt}{\bar{R}_j(t)}
\newcommand{\xit}{x_i(t)}
\newcommand{\xjt}{x_j(t)}
\newcommand{\uxit}{\bar{x}_i(t)}
\newcommand{\ti}{\theta_i}
\newcommand{\tit}{\theta_i(t)}
\newcommand{\tjt}{\theta_j(t)}
\newcommand{\htit}{\hat{\theta}_i(t)}
\renewcommand{\tt}{\bm{\theta}(t)}
\newcommand{\htt}{\hat{\bm{\theta}}(t)}
\newcommand{\lit}{l_i(t,\bm{\theta}(t))}
\newcommand{\hlit}{l_i(t,\htt)}
\newcommand{\xkt}{\bm{x}_k(t)}
\newcommand{\xt}{\bm{x}(t)}
\newcommand{\E}{\mathbb{E}}
\newcommand{\pkt}{\mathsf{P}_k(t)}
\newcommand{\hpkt}{\hat{\mathsf{P}}_k(t,\htt)}
\newcommand{\fg}{f_{ij}^{\text{ground}}}
\newcommand{\fa}{f_{ij}^{\text{air}}}

\newtheorem{prop}{Proposition}

\allowdisplaybreaks

\begin{document}

\begin{frontmatter}
\title{Vaccine allocation policy optimization and budget sharing mechanism using Thompson sampling}

\author[SKEMA]{David Rey}
\author[UNSW2]{Ahmed W. Hammad}
\author[UNSW]{Meead Saberi\corref{mycorrespondingauthor}}
\ead{meead.saberi@unsw.edu.au}
\cortext[mycorrespondingauthor]{Corresponding author}
\address[SKEMA]{SKEMA Business School, Universit\'e C\^ote d'Azur, Sophia Antipolis, France.}
\address[UNSW2]{School of Built Environment, UNSW Sydney, Sydney, NSW, 2052, Australia}
\address[UNSW]{School of Civil and Environmental Engineering, UNSW Sydney, Sydney, NSW, 2052, Australia}

\begin{abstract}
The optimal allocation of vaccines to population subgroups over time is a challenging health care management problem. In the context of a pandemic, the interaction between vaccination policies adopted by multiple agents and the cooperation (or lack thereof) creates a complex environment that affects the global transmission dynamics of the disease. In this study, we take the perspective of decision-making agents that aim to minimize the size of their susceptible populations and must allocate vaccine under limited supply. We assume that vaccine efficiency rates are unknown to agents and we propose an optimization policy based on Thompson sampling to learn mean vaccine efficiency rates over time. Furthermore, we develop a budget-balanced resource sharing mechanism to promote cooperation among agents. We apply the proposed framework to the COVID-19 pandemic. We use a raster model of the world where agents represent the main countries worldwide and interact in a global mobility network to generate multiple problem instances. Our numerical results show that the proposed vaccine allocation policy achieves a larger reduction in the number of susceptible individuals, infections and deaths globally compared to a population-based policy. In addition, we show that, under a fixed global vaccine allocation budget, most countries can reduce their national number of infections and deaths by sharing their budget with countries with which they have a relatively high mobility exchange. The proposed framework can be used to improve policy-making in health care management by national and global health authorities. 
\end{abstract}


\end{frontmatter}

\section{Introduction}
\label{intro}

With the continued spread of the coronavirus disease 2019 (COVID-19) around the globe and an increasing number of administered vaccines, discussions on the safe, effective, and ethical allocation of COVID-19 vaccines are growing \citep{Emanuel2020,LIBOTTE2020105664,Persad2020,Peiris2020,Bollyky2020,NAP25917}. The global allocation and distribution of COVID-19 vaccines is a challenging logistical problem \citep{Roope2020,ZAFFRAN2013B73}. The effect of human mobility on the spread of disease is also well known \citep{kraemer2020effect}. The need for a fair and ethical allocation of vaccines further reinforces this challenge given the emerging controversial issues on public health, diplomacy and economics \citep{Liu499}. Furthermore, the interaction between the vaccination policies adopted by different countries and vaccine protectionism creates a complex environment that affects the global transmission dynamics of COVID-19 and the efficiency of vaccination strategies. Furthermore, vaccine protectionism among a few developed countries has slowed down the global vaccination efforts against COVID-19 \citep{nkengasong2021covid,lancet2020global}. 

The inclination of some of the vaccine producer countries toward vaccine protectionism has negatively affected the global progress against COVID-19. This will likely limit access of low-income countries to vaccines \citep{Amnesty}. In response, the World Health Organization (WHO) has initiated a global collaborative initiative known as COVAX \citep{WHO} as part of the Access to COVID-19 Tools (ACT) Accelerator that aims to accelerate the equitable access to COVID-19 vaccines. Despite a slow rollover, COVAX aims to offer vaccines to countries in amounts proportional to their population starting with 3\% of each country’s population and gradually increasing the allocation to at least 20\%. While a population-based vaccine allocation policy may appear equitable, there are inherent limitations to this approach that overlooks the global transmission dynamics of SARS-CoV-2, worldwide and local human mobility, and countries’ vastly different vaccination capacities. Further, the true efficiency of COVID-19 vaccines may not match their efficiency observed under clinical conditions. The timely question of what the optimal allocation strategy of COVID-19 vaccines at the global scale is remains open. The true cost of vaccine protectionism is also not yet well-understood.

In this study, we examine the problem of allocating vaccines across populations in the presence of uncertainty with regards to the efficiency of vaccination decisions. We take the perspective of a decision-making agent (e.g. country) that is in charge of vaccinating its population. We assume that the agent's population is spread over space and that, given a vaccination budget, the agent must decide how to allocate vaccines across its population. We consider discrete time and we assume that agents make vaccine allocation decisions periodically (e.g. weekly), based on the distribution of infected individuals across their population. We assume that the impact of vaccine allocation decisions is unknown to agents and that agents learn mean vaccine efficiency rates by observing the impact of past decisions. We frame this complex problem as an online learning and optimization problem. We propose a data-driven optimization method based on Thompson sampling (TS) -- a reinforcement learning approach that uses Bayesian optimization to learn vaccine efficiency rates over time -- to identify competitive vaccine allocation policies. We also propose a budget-balanced resource sharing mechanism that aims to further mitigate the impact of global mobility patterns and promote resource sharing across agents.

To test the proposed approach, we conduct a global numerical study wherein countries of the world are represented as agents. To capture disease dynamics, we embed a metapopulation epidemic model with a raster model of the world where each cell is a node in a global mobility network. We implement the proposed TS-based vaccine allocation policy along with several benchmarking policies over a two-year horizon. Our findings provide evidence that a population-based vaccine allocation policy is sub-optimal compared to the proposed TS-based policy. We show that the latter achieves a larger reduction in the number of infected and death globally. We also provide supporting evidence that a budget sharing mechanism between agents could further reduce the number of infections and deaths globally. Under fixed vaccine allocation budget, we show that most countries can reduce their national infection and mortality rates by sharing their budget with countries with which they have a relatively high mobility exchange. This counter-intuitive finding contradicts the popular belief on the national (as opposed to international) benefits of vaccine protectionism and reveals the significant potential impact of cooperation among countries.

The rest of the paper is organized as follows. We review the literature on vaccine allocation in Section \ref{lit}. We present the global metapopulation epidemic model used to represent disease spreading in Section \ref{pre}. The proposed online vaccine allocation problem along with the TS-based policy are introduced in Section \ref{vap}. The budget sharing mechanism is presented in Section \ref{sharing}. The data of the global case study is summarized in Section \ref{data}, and the numerical results are provided in Section \ref{num}. Concluding remarks are discussed in Section \ref{con}. 

\section{Literature Review}
\label{lit}

We review the state-of-the-art in the field of vaccine allocation. We start by examining studies which have discussed the modeling of vaccine allocation and vaccination practice in Section \ref{vamod}. We then focus on studies that have proposed mathematical optimization formulations to solve vaccine allocation problems in Section \ref{vaopt} before outlining our contributions with respect to the literature in Section \ref{contributions}.

\subsection{Vaccine allocation modeling and practice}
\label{vamod}

Numerous studies in the literature have explored different vaccination strategies under limited resources including prioritizing the allocation based on age \citep{Emanuel2020,Couzin960,Lipsitch41}, demographics \citep{BECKER1997117}, underlying high-risk conditions \citep{Tuite2010}, and virus transmission dynamics \citep{Medlock1705}. However, very few studies used a computational model of global epidemic and human mobility. One of the more widely known computational models used to study the transmission of infectious diseases is the Global Epidemic and Mobility (GLEaM) model \citep{van2011gleamviz} which is a stochastic computational model that integrates worldwide human mobility data and high-resolution demographic to simulate spread of disease at the global scale. For example, GLEaM has been previously used to model and assess the effectiveness of different H1N1 vaccination campaigns \citep{Bajardiet2009}.

The allocation of vaccines to countries proportionally to their population could result in significantly different morbidity and mortality rates in equally populous countries \citep{Emanuel2020} because of varying transmission dynamics of the virus due to differences in age structure, health care capacity, logistics infrastructure, and implemented non-pharmaceutical interventions in each country. Thus, finding the optimal global allocation of COVID-19 vaccines requires a model that simultaneously accounts for the virus transmission dynamics, human mobility, vaccines’ efficiency data, vaccination capacity, and more importantly the interaction between the vaccine allocation strategies adopted by different countries and their impact on the number of infections and deaths globally. 

Many existing studies on vaccine allocation are limited to sensitivity analyses \citep{Tuite2010,MYLIUS20083742} with very few formulating vaccine allocation problems using mathematical optimization and at a global scale. In addition, the vast majority of existing vaccine allocation strategies are static and established over a pre-defined time horizon, i.e., they do not allow the possibility to dynamically adapt vaccine allocations based on new data. 

\subsection{Optimizing vaccine allocation policies}
\label{vaopt}

Research on the optimization of vaccine allocation policies appear to have emerged in the field of Operations Research (OR) with the study of \citet{longini1978optimization} who presented a vaccine allocation formulation to minimize the total number of infections under limited supply. \citet{BECKER1997117} proposed a linear programming formulation to identify the optimal allocation of vaccines in a community of households. In their study, the authors seek to either minimize vaccination coverage or the reproduction number while accounting for disease transmission dynamics. \citet{ball2002optimal} extend this research by considering two level of population mixing, i.e. within households and among the population of households and discuss optimality conditions for vaccine allocation policies in this context. \citet{HILL200385} proposed optimal vaccination policies for populations divided into subgroups. \citet{TANNER2008144} build on prior research on epidemics within population of households by modeling the upper bound on the reproduction number using chance-constraints. This framework is further refined by \citet{tanner2010iis} who developed a customized branch-and-cut algorithm to improve computational scalability. This stream of research has mainly focused on allocating vaccines among subgroups of a population which are often age-stratified. More recently, \citet{enayati2020optimal} built on this stream of research and developed mathematical optimization formulation to incorporate equity within the vaccine allocation process. In their study, the authors consider population subgroups based on both location and age.

In the context of the H1N1 pandemic, \citet{samii2012reservation} explored the problem of determining optimal reservation and allocation policies for vaccine inventory rationing. Several researchers have also studied the problem of allocating vaccine among subgroups of a population where subgroups represent geographical regions. Unlike previous studies, \citet{teytelman2013multiregional} addressed a dynamic variant of the vaccine allocation problem and the authors developed heuristics to find optimal policies in this context. \citet{YARMAND2014208} proposed a two-stage stochastic programming approach for vaccine allocation over multiple locations where the second stage can be viewed as a recourse stage which is triggered if the outbreak is not contained. \citet{long2018spatial} adopts a multi-period stochastic optimization framework and aims to identify the optimal allocation of health care resources across space. 

The COVID-19 pandemic has reinvigorated research on vaccine allocation problems, especially in the context of spatial and temporal decisions. \citet{yang2021optimizing} addresses the problem of optimizing vaccine distribution networks in low and middle income countries wherein the lack of health care resources can substantially affect the efficiency of the vaccine distribution chain \citep{de2020vaccine}. \citet{bertsimas2020optimizing} used a DELPHI compartmental model to capture disease dynamics and proposed a simulation-based optimization approach to solve this vaccine allocation problem. A similar modeling framework was also used by \citet{bertsimas2021locate} to identify the optimal location of COVID-19 mass vaccination facilities at the country (US) level. \citet{chen2020allocation} propose both static and dynamic vaccine allocation policies in an age-stratified population and focus on a case study based New York City data. \citet{thul2021stochastic} proposed a stochastic optimization for vaccine and testing kit allocation. The authors adopt an online learning and optimization context where the decision-maker must repeatedly decide how to allocate vaccines and testing kits across space. They model this problem as a partially-observable Markov Decision Process and develop several policies. Their computational study is focused on the US wherein states represent population subgroups. 

\newpage

\subsection{Our contributions}
\label{contributions}

While the COVID-19 pandemic has triggered several new studies on OR-driven vaccine allocation methodologies there remain significant research gaps in the literature. The vast majority of studies have focused on static vaccine allocation problems, i.e. the decision-maker seeks to find the optimal policy to allocate vaccines -- possibly over multiple time periods -- across a population based on the available information at the time of decision. Only a few studies considered either recourse actions using a stochastic programming approach \citep{teytelman2013multiregional,YARMAND2014208,chen2020allocation}, or an online resource allocation framework wherein decisions needs to be taken repeatedly over time in light of new data. Among these, the study of \citet{thul2021stochastic} is the only one which consider a learning and optimization framework to iteratively refine the vaccine allocation policy based on historical observations. Our study adopts a related framework but differs in the choice of the uncertain parameters modeled and in the proposed learning and optimization methods. Furthermore, most studies have only examined the vaccine allocation problem from the perspective of a single decision-making agent, e.g. a single country or a local health authority. Although comprehensive global case studies have been conducted in the field of human mobility and epidemic modeling \citep{Liu499,Bollyky2020}, the interaction among several agents at a global scale has remained largely unaddressed in the health care management literature on vaccine allocation problems. 

In this study, we attempt to address some of these research gaps. We consider the problem of allocating vaccines to a population of subgroups where the latter represent geographical regions while accounting for global mobility and disease dynamics. Unlike previous studies, we do not restrain our analysis to a single decision-maker. Instead, we consider multiple decision-making agents and each agent aims to allocate vaccines to minimize the size of its susceptible population. We believe that this multi-agent configuration is representative of ongoing efforts to mitigate the COVID-19 pandemic by countries worldwide. We adopt an online optimization framework wherein each agent is assumed to periodically make vaccine allocation decisions subject to a budget constraint based on historical data. To capture the interaction among agents, we construct a global mobility network and embed a compartmental epidemic model to represent disease dynamics within each population subgroup. We assume that vaccine allocation decision only makes a proportion of the susceptible population immune. We refer to this proportion as the vaccine efficiency rate. We assume that vaccine efficiency rates are unknown to agents and that the latter learn these rates over time based on past decisions. We develop a customized reinforcement learning policy based on TS to dynamically solve agent-level vaccine allocation problems. 

To test the proposed vaccine allocation policies, we generate a raster model of the world along with mobility models to obtain global vaccine allocation problem instances where agents represent countries worldwide. We also propose a budget sharing mechanism to explore the potential benefits of cooperation among agents. While the popular belief is that vaccine protectionism is in the interest of vaccine producer countries, in this study, we provide evidence that most countries could further reduce their national death toll by sharing their (fixed) vaccine allocation budget with connected countries via the global human mobility network. This counter-intuitive finding demotes vaccine protectionism and promotes a more equitable allocation strategy across the world in which countries with higher GDP per capita can in fact benefit from sharing their fixed vaccine allocation budget with lower GDP per capita countries with which they have a higher mobility exchange. 

\section{Preliminaries}
\label{pre}

We next present our modeling approach in two parts. We first introduce the global mobility network in Section \ref{mob} before describing the metapopulation epidemic model in Section \ref{epi}.

\subsection{Global mobility network model}
\label{mob}

To model the global mobility patterns, we use a uniform raster model of the world wherein each population cell is represented by a node in a network. Let $\V$ be the set of nodes in this network. Each population in the model is represented by a node $i \in \V$. Let $\Ni \subset \V$ be the neighborhood of node $i$ which represents the nodes connected to $i$ in the network. Arcs among pairs of nodes in the network are introduced for all pairs with non-zero mobility flows. We assume that individuals can move between nodes via ground and/or air connections. We next describe how ground and air mobility flows are determined and combined into a global flow matrix. 

We use the radiation model \citep{simini2012universal} to determine ground mobility flows between nodes of the network. Let $d_{ij}$ representing the great circle distance between nodes $i,j \in \V$, and let $D$ be a distance threshold beyond which we assume that ground mobility is null. We define the ground neighborhood of node $i \in \V$ as $\NGi \equiv \{j\in\V: d_{ij}\leq D\}$, for all $i \in \V$. Let $P_i$ be the population of node $i \in \V$ and let $\alpha_i$ be the fraction of population at node $i$ that commutes. The ground mobility flow from node $i \in \V$ to $j \in \NGi$, denoted $\fg$, can be determined using Eq. \eqref{eq:ground}. 

\begin{equation}
\fg \equiv \alpha_i P_i \frac{P_i P_j}{\sum\limits_{i' \in \NGi} P_{i'} \left(P_j + \sum\limits_{j' \in \NGi} P_{j'}\right)}, \quad \forall i \in \V, j \in \NGi.
\label{eq:ground}
\end{equation}

To model air mobility, we generate a Voronoi polygon around each international airport in the world and identify the set of nodes that fall inside each Voronoi polygon. We then proportionally distribute the inflows to and outflows from each airport to its corresponding set of nodes based on node populations. Let $\C$ be the set of airports also referred to as the set of Voronoi polygons. Formally, the Voronoi polygon corresponding to airport $a$ denoted $\Pi(a)$ is defined as a subset of $\V$ as:
$\Pi(a) \equiv \{i \in \V : d_E (i,a)\le d_E (i,b), \quad \forall a,b\in\C, a\neq b\}$, where $d_E$ is a distance function that represents the Euclidean distance on $\mathbb{R}^2$. For each node $i \in \V$, we denote $\mu_i \in \C$ its assigned airport. Given a matrix of air flows among worldwide airports $[g_{ab}]_{a,b \in \C}$, we determine the air neighborhood of node $i \in \V$ as the set of nodes $j \in \V$ that are connected to $i$ via an air link. Let $\NAi \equiv \{j \in \V : g_{\mu_i \mu_j} > 0\}$, for all $i \in \V$ be the air neighborhood of node $i$ and $P_a = \sum_{i \in \Pi(a)} P_i$ be the population of the Voronoi polygon $a\in \C$. We determine the air mobility flow from node $i$ to $j$, denoted $\fa$, as:
\begin{equation}
    \fa \equiv g_{\mu_i \mu_j} \frac{P_i + P_j}{P_{\mu_i} + P_{\mu_j}}, \quad \forall i \in \V, j \in \NAi.
\end{equation}

We define $\Ni \equiv \NGi \cup \NAi$ as the global neighborhood of node $i \in \V$ and $f_{ij} \equiv \fg + \fa$ as the mobility flow between two nodes $i,j \in \V$. We refer to $[f_{ij}]_{i\in \V, j \in \Ni}$ as the global flow matrix.

\subsection{Metapopulation epidemic model}
\label{epi}

We model the spread of diseases using a metapopulation discrete-time compartmental epidemic model. We adopt the generic metapopulation model of \citet{brockmann2013hidden} and use a Susceptible-Infected-Recovered-Dead (SIRD) model to estimate the global spreading dynamics of the disease. This model is an extension of the classical Susceptible-Infected-Recovered (SIR) model \citep{kermack1927contribution} where the dead compartment (D) represents the fraction of infected (I) individual which are expected to die from the disease.

We define a flow rate matrix $[\pij]_{i\in \V, j\in \Ni}$ based on the flow matrix $[\fij]_{i\in \V, j\in \Ni}$ as follows:
\begin{equation}
\pij \equiv \frac{\fij}{\sum\limits_{j' \in \Ni} f_{ij'}}, \quad \forall i \in \V, j \in \Ni.
\end{equation}

In addition, we define the global flow-to-population ratio $\rho$ as:
\begin{equation}
\rho \equiv \frac{\sum\limits_{i \in \V} \sum\limits_{j \in \Ni} \fij}{\sum\limits_{i\in \V} P_i}.
\end{equation}

Local disease transmission rate $\beta_i$, recovery rate $\gamma_i$, and death rate $\lambda_i$ for each node $i \in \V$ are assumed to be known. Let $P_i$ be the population of node $i \in \V$. For any compartment $U = S$, $I$, $R$ or $D$, we define compartment proportions $\bar{U}_i(t)$ as $\bar{U}_i(t) \equiv U_i(t)/P_i$. Let $\T$ be the set of time periods. Assuming a constant population (including deaths), using the flow rate matrix $[\pij]_{i\in \V, j\in \Ni}$ and the global flow-to-population $\rho$, the SIRD model with mobility at node $i \in \V$ can be represented by the system of equations:
\begin{subequations}
\begin{align}
&\hsitp = \hsit - \beta_i \hsit\hiit + \rho \sum_{j \in \Ni} \pij (\hsjt - \hsit), &&\forall t \in \T,\\
&\hiitp = \hiit + \beta_i \hsit\hiit - \gamma_i\hiit + \rho \sum_{j \in \Ni} \pij (\hijt - \hiit), &&\forall t \in \T,\\
&\hritp = \hrit + (1 - \lambda_i)\gamma_i \hiit + \rho \sum_{j \in \Ni} \pij (\hrjt - \hrit), &&\forall t \in \T, \\
&\hditp = 1 - \hsitp - \hiitp - \hritp, &&\forall t \in \T.
\end{align}
\label{eq:SIRD}
\end{subequations}

The metapopulation model Eq. \eqref{eq:SIRD} is based on two underlying assumptions: i) for any pair of nodes $i,j \in \V$, $\pij P_i = \pji P_j$ and ii) for any node $i \in \V$, outflow is proportional to population: $P_i \sim \sum_{j \in \Ni} \fij$ as discussed by \citet{brockmann2013hidden}. We discuss to which degree these assumptions are verified by the data of this case study in Section \ref{data}.\\

We next introduce the proposed online vaccine allocation problem, mathematical optimization formulations and online learning algorithms.

\section{Online vaccine allocation problem}
\label{vap}

We model the vaccination of populations as an online optimization problem. We consider a set of decision-making agents that must allocate vaccines to their population over a set of time periods. (In our case study, agents represent countries.) Each agent controls a subset of nodes of the network. At each time period, agents are assumed to have limited resources for allocating vaccines and must decide how to allocate these resources across their populations. We assume that allocating vaccines at a node makes immune a proportion of the susceptible population of this node and we call this proportion the vaccine efficiency rate. We model the impact of vaccine allocation decisions as a stochastic process: we assume that vaccine efficiency rates are unknown to agents and that these rates are learned over time when observing the impact of vaccine allocation decisions.

We first define vaccine allocation decisions and introduce the proposed stochastic optimization framework within the metapopulation epidemic model in Section \ref{stoch}. We then present online resource allocation problem faced by agents along with a mathematical optimization formulation in Section \ref{form}. We propose an online learning algorithm based on Thompson sampling to identify optimal policies for vaccine allocation in Section \ref{algo}. 

\subsection{Modeling vaccination allocation decisions}
\label{stoch}

Let $\xit \in [0,1]$ be a real decision variable representing the proportion of vaccines allocated at node $i \in \V$ at time $t \in \T$, i.e. $\xit=1$ corresponds to the case where node $i$ is fully supplied in vaccines whereas $0 < \xit < 1$ corresponds to a partial allocation, and $\xit=0$ means that no vaccines at allocated at $i$ at time $t$. Let $\ti \in [0,1]$ be a random variable representing the mean vaccine efficiency rate of node $i \in \V$ and let $\bm{\theta}$ be the vector of mean vaccine efficiency rates. We assume that $\bm{\theta}$ and the probability distributions of vaccine efficiency rates are unknown to agents and are observed after making vaccine allocation decisions. Let $\tit$ be the observed vaccine efficiency rate at node $i \in \V$ at time $t \in \T$. When designing vaccination allocation policies, for each node $i \in \V$, we will later require that random variables $\tit$ follow arbitrary probability distributions with support in $[0,1]$ and mean $\ti$. To capture the impact of vaccine allocation decisions, we assume that the susceptible population $\sit$ is reduced to $(1 - \tit\xit)\sit$ and that $\tit\xit\sit$ individuals are moved to the recovered compartment. Hence, we incorporate vaccine allocation decisions within the proposed SIRD model \eqref{eq:SIRD} as follows:
\begin{subequations}
\begin{align}
&\hsitp = \left(\hsit - \beta_i \hsit\hiit\right)\left(1 - \tit\xit\right) \nonumber \\
& \qquad\qquad\quad + \rho \sum_{j \in \Ni} \pij \left((\hsjt(1 - \tjt\xjt) - \hsit(1 - \tit\xit)\right), &&\forall t \in \T, \label{eq:Sx}\\
&\hiitp = \hiit + \beta_i \hsit\hiit \left(1 - \tit\xit\right) - \gamma_i\hiit + \rho \sum_{j \in \Ni} \pij (\hijt - \hiit), &&\forall t \in \T,\\
&\hritp = \hrit + \hsit\tit\xit + (1 - \lambda_i)\gamma_i \hiit \nonumber \\
& \qquad\qquad\quad + \rho \sum_{j \in \Ni} \pij \left(\hrjt + \hsjt\tjt\xjt - \hrit - \hsit\tit\xit\right), &&\forall t \in \T, \\
&\hditp = 1 - \hsitp - \hiitp - \hritp, &&\forall t \in \T.
\end{align}
\label{eq:SIRDx}
\end{subequations}

\subsection{Vaccine allocation formulation}
\label{form}

Let $\K$ be the set of decision-making agents. Each agent $k \in \K$ controls a set of nodes $\Vk \subset \V$. We assume that the goal of agents is to minimize the cumulative number of susceptible people by solving a sequence of online vaccine allocation problems -- one per time period -- over a given time horizon. We next describe the vaccine allocation problem solved at time $t \in \T$ by agent $k \in \K$.\\

We propose a direct lookahead approximation (DLA) policy to minimize the expected number of susceptible individuals at the next time period \emph{(29)}. We assume that global data for time period $t$, including compartment volumes $\hsit$, $\hiit$, $\hrit$ and $\hdit$, are available to agents when allocating vaccines for time period $t+1$. The decisions of other agents are unknown to agent $k$. Given an agent $k \in \K$, let $-k$ denote the other agents and let $\Vnk$ represent the nodes controlled by other agents in the network. We take a worst-case approach to capture this competitive effect and thus assume that the nodes $j \in \Vnk$ are not allocated any vaccines i.e. $\xjt=0$. Let $\xkt \in [0,1]^{|\Vk|}$ be the vector of decision variables of agent $k \in \Vk$ at time $t \in \T$. The objective function of agent $k$ at time $t$ is to minimize the expected number of susceptible individuals at time $t+1$:
\begin{align}
\text{minimize}\ &\E\Bigg\{\sum_{i \in \Vk} \Bigg(\left(\hsit - \beta_i \hsit\hiit\right)\left(1 - \tit\xit\right) \nonumber\\
& \quad + \rho \sum_{j \in \Ni \cap \Vk} \pij \left((\hsjt(1 - \tjt\xjt) - \hsit(1 -
\tit\xit)\right) \nonumber \\
& \quad + \rho \sum_{j \in \Ni \cap \Vnk} \pij \left((\hsjt - \hsit(1 - \tit\xit)\right)\Bigg) \Bigg\}. \label{eq:obj}
\end{align}

The expectation in Eq. \eqref{eq:obj} is taken over the random variables $\tit$ for all $i \in \Vk$. From a decision-making standpoint, the number of susceptible individuals can be viewed as potential loss and the goal of Eq. \eqref{eq:obj} is to minimize the expected loss. For conciseness, it is convenient to rewrite Eq. \eqref{eq:obj} in compact form by eliminating constants and aggregating all the coefficients of variable $\xit$ in a weight $\lit$ representing the loss of node $i$ at time $t$. Observe that $\lit$ depends on random variables $\tt$ and thus is also a random variable. Accordingly, the objective function Eq. \eqref{eq:obj} is rewritten compactly as:
\begin{equation}
\text{minimize}\ \E\Bigg\{\sum_{i \in \Vk} \lit \xit \Bigg\}.
\label{eq:obj2}
\end{equation}

Let $\Gamma_k$ be the per-period vaccination capacity of agent $k$. This capacity represents the ability of an agent to distribute vaccines to its population at each time period. Let $B_k(t)$ be the budget of agent $k \in \K$ at time period $t \in \T$ for allocating vaccines across the set of nodes $\Vk$. We assume that agents' per-period budget is function of their vaccination capacity, i.e. $B_k(t) = f(\Gamma_k)$. The per-unit cost to allocate vaccines at node $i \in \Vk$ is assumed to require a known cost $C_i$. Data used to generate values for $\Gamma_k$, $B_k(t)$ and $C_i$ are discussed in Section \ref{capdata}. At time period $t \in \T$, the budget constraint of agent $k$ is:
\begin{equation}
\sum_{i \in \Vk} \xit C_i \leq B_k(t).
\end{equation}

Even though node populations (including the death compartment) are assumed to be constants, global mobility across the network means that the number of susceptible individuals at a node might fluctuate in such a way that after a certain number of time periods the susceptible population of a node might be composed of individuals which were not located at this node at the beginning of the time horizon under consideration. This implies that vaccine allocation decisions across space might need to be repeated. This is particularly critical if the ratio of the inflow to the population of a node is relatively large. Alternatively, allocating vaccines at nodes at which inflow is low relative to population may not require additional vaccines before several time periods. Hence, it is critical to account for the history of vaccine allocation decisions in the proposed modeling framework. To capture this effect, we assume that at each time period $t$, the upper bound  $\uxit \leq 1$ on $\xit$ is determined based on the history of vaccine allocation decisions over $t' \in \{0,\ldots,t-1\}$. The update rule to determine $\uxit$ is discussed in Section \ref{algo} which outlines the proposed TS-based algorithm for the online allocation of vaccines. The vaccine allocation problem of agent $k \in \K$ at time $t \in \T$ is denoted $\pkt$:
\begin{subequations}
\begin{align}
\pkt: \quad & \text{minimize}\  \E\Bigg\{\sum_{i \in \Vk} \lit \xit \Bigg\}, \\
& \text{subject to:} \nonumber\\
& \sum_{i \in \Vk} \xit C_i \leq B_k(t), \\
& 0 \leq \xit \leq \uxit, \forall i \in \Vk.
\end{align}
\label{eq:pkt}
\end{subequations}

The optimization problem $\pkt$ can be viewed as an online knapsack problem where the ``cost'' (loss) of items (nodes) are stochastic and depend on unknown vaccine efficiency rates. 

\subsection{Online vaccine allocation by Thompson sampling}
\label{algo}

To solve $\pkt$ we propose a reinforcement learning approach based on Bayesian optimization. We adapt the algorithm proposed by \citet{thompson1933likelihood}, also known as Thompson Sampling (TS), to account for resource allocation constraints. TS has been shown to be an efficient algorithm for DLA policies and empirical studies have shown that TS is highly competitive to address the exploration-exploitation tradeoff in online learning problems \citep{chapelle2011empirical}. TS has also been adapted to constrained online optimization problems such as linear-quadratic control \citep{abeille2018improved}, online network revenue management \citep{ferreira2018online} and real-time energy pricing \citep{tucker2020constrained}.\\

We next adapt the TS algorithm proposed by \citet{agrawal2012analysis} for general stochastic bandits to solve the online vaccine allocation problem at hand. This TS algorithm only requires to assume that mean vaccine efficiency rates $\bm{\theta}$ are generated from an arbitrary unknown distribution with support in $[0,1]$ which fits well the purpose of this study, i.e. learning mean vaccine efficiency rates. In addition, this TS algorithm uses Beta distributions as Bayesian priors and we adopt the same framework to model agents' beliefs over mean vaccine efficiency rates. Accordingly, for each node $i \in \V$, we denote $\mathsf{Beta}(a_i,b_i)$ the prior of its mean vaccine efficiency rate where $a_i$ and $b_i$ are parameters of the Beta distribution. Let $\htt$ denote the vector of mean vaccine efficiency rates sampled from priors $\mathsf{Beta}(a_i,b_i)$ for all $i \in \V$ at time $t \in \T$. Given time period $t \in \T$, we denote $\htit$ the mean vaccine efficiency rate of node $i \in \V$ sampled from prior $\mathsf{Beta}(a_i,b_i)$, and we denote $\hlit$ the corresponding sampled loss function, i.e. $\hlit$ is the coefficient of variable $\xit$ in Eq. \eqref{eq:obj} obtained by substituting random variables $[\tit]_{i\in\V}$ with sampled mean vaccine efficiency rates $[\htit]_{i\in \V}$. The approximated vaccine allocation problem of agent $k \in \K$ at time $t \in \T$ is denoted $\hpkt$:
\begin{subequations}
\begin{align}
\hpkt: \quad & \text{minimize}\  \sum_{i \in \Vk} \hlit \xit, \label{eq:hobj}\\
& \text{subject to:} \nonumber\\
& \sum_{i \in \Vk} \xit C_i \leq B_k(t), \\
& 0 \leq \xit \leq \uxit, \forall i \in \Vk.
\end{align}
\label{eq:hpkt}
\end{subequations}

Formulation $\hpkt$ is a linear knapsack problem that can be solved in polynomial-time using a greedy algorithm by sorting nodes $\Vk$ by increasing loss-to-cost ratio $[\hlit/C_i]_{i \in \Vk}$.\\

The pseudo-code of the proposed TS-based vaccine allocation policy is summarized in Algorithm \ref{alg:TS}. In our numerical experiments, all parameters $a_i$ and $b_i$ of the prior distribution of nodes $i \in \V$ are initialized to 1 which corresponds to uniform distributions (lines 2 and 3). In practice, historical vaccine allocation data could be used to improve the initialization of prior distribution parameters. At each time period $t$, prior distributions $\mathsf{Beta}(a_i,b_i)$ are sampled to obtain estimates of mean vaccine efficiency rates $[\htit]_{i\in\V}$ (line 6). Then, for each agent $k\in\K$, the approximated vaccine allocation problem $\hpkt$ is solved to determine the vaccine allocation strategy $\xkt$ (line 8). For each node $i\in\V$, the upper bound $\uxit$ is then updated by deducing the amount of vaccines allocated to this node over the time window $[\max\{t-m_i+1,1\},t]$ (line 10), where $m_i$ is a node-based parameter used to adjust the width of historical observations taken into consideration at each decision epoch. In our numerical experiments, $m_i$ is determined as the ratio of the node population $P_i$ to the total inflow at this node $\sum_{j \in \Ni} f_{ji}$ rounded to the nearest integer above. Hence, the parameter $m_i$ represents a conservative estimate of the number of time periods needed to ``renew" the population from a human mobility standpoint, and the time window $[\max\{t-m_i+1,1\},t]$ represents the history of allocation decisions taken into consideration for updating the upper bounds $\uxit$. Random variables representing vaccine efficiency rates $[\tit]_{i\in\V : \xit >0}$ are observed for all nodes that are allocated vaccines, and the vaccine-dependent SIRD model represented by \eqref{eq:SIRDx} is then solved to obtain node compartments for the next time period (line 11). Bernoulli trials using the observed vaccine efficiency rates are performed (line 14) and their outcomes are used to update the parameters of the corresponding prior distributions on vaccine efficiency rates (line 14-18) for the next time period.\\

\begin{algorithm}[ht]
\label{alg:TS}
\For{$i \in \V$}{
    $a_i \gets 1$\\
    $b_i \gets 1$\\
}
\For{$t \in \T$}{
    \For{$i \in \V$}{
        $\htit \gets$ sample prior $\mathsf{Beta}(a_i,b_i)$\\
    }
    \For{$k \in \K$}{
        $\xkt \gets$ solve $\hpkt$\\
        \For{$i \in \Vk$}{
            $\uxit \gets 1 - \sum_{t' = \max\{t - m_i + 1,1\}}^{t} x_i(t')$\\
        }
    }
    $[\hsitp,\hiitp,\hritp,\hditp]_{i
    \in \V} \gets$ solve \eqref{eq:SIRDx} using $\xt$ and observe $[\tit]_{i\in\V : \xit >0}$ \\
    \For{$i \in \V$}{
        \If{$\xit > 0$}{
            $y_i \gets$ Bernoulli trial $(\tit)$\\
            \If{$y_i = 1$}{
                $a_i \gets a_i + 1$ \\
            }
            \Else{
                $b_i \gets b_i + 1$ \\
            }
        }
    }
}
\caption{TS-based policy for online vaccine allocation}
\end{algorithm}

\subsection{Benchmarking with other vaccine allocation policies}
\label{bench}

We compare the proposed TS-based vaccine allocation policy (TS) with three alternative approaches: a population-based (PB) approach, a moving average (MA) approach, and a greedy learning (GY) approach. These three methods are described below and their difference with the proposed TS-based approach are discussed.

\begin{itemize}
    \item PB allocates vaccines using node population sizes to measure the expected impact of allocation decisions. This is equivalent to replace the loss function $l_i(t,\hat{\bm{\theta}}(t))$ of node $i$ at time $t$ by $P_i$ in the objective function \eqref{eq:hobj}. Since $P_i$ does not depend on historical data, this strategy is the easiest to implement as it does not require any tracking or learning of stochastic parameters.
    
    \item MA allocates vaccines by estimating the loss function \eqref{eq:hobj} $l_i(t,\hat{\bm{\theta}}(t))$ of node $i$ at time $t$ using a moving average over the historical data. Hence, instead of using a learning approach to estimate unknown vaccine efficiency rates as in the TS-based approach, the MA approach simply estimates node-based vaccine efficiency rates as the average of the observed data at this node.
    
    \item GY allocates vaccines using a Bayesian optimization approach, similarly to the proposed TS-based approach. The only difference between GY and TS is that the former uses the expected value of the prior distribution instead of sampling from the prior distribution when estimating nodes' mean vaccine efficiency rates (line 6). 
\end{itemize}

\section{Budget-balanced resource sharing mechanism}
\label{sharing}

The proposed vaccine allocation policies outlined in Section \ref{vap} assume that all agents act independently, without sharing any vaccination resources. Here, we propose a budget sharing mechanism that aims to further improve the impact of vaccine allocation decisions by re-distributing vaccination budgets across agents. We show that the proposed mechanism is budget-balanced, i.e. there is no additional budget incorporated in the system, hence this sharing mechanism can be compared to a ``no sharing'' mechanism to measure its efficiency. The proposed budget-balanced resource sharing mechanism tracks, for each agent, the ratio of internal versus external infections. At every time period, we determine the proportion of budget shared with other agents at the next time period as the proportion of external infections to total, i.e., internal and external, infections. The shared budget is then split among connected agents via the global mobility network proportionally to their volume of external infections weighted by agents’ vaccination capacities. Hence, the proposed mechanism promotes a more equitable allocation of vaccination resources across agents by ensuring that agents with a relatively low vaccination capacity receive a proportionally larger share of budget compared to agents with a high vaccination capacity.

To implement the proposed budget sharing mechanism, at the beginning of each time period $t$, we update agent-based vaccine allocation budgets $B_k(t)$ by determining the amount of budget shared proportionally to the ratio of external to total, i.e. external and internal, infections. We then allocate the portion of budget shared to connected agents in the global mobility network proportionally to the volume of infected population weighted by agents' vaccination capacities.\\

Formally, at time period $t$, for each agent $k\in\K$ and for each node $i\in \Vk$, let $\bar{I}^{\text{in}}_i(t+1)$ and $\bar{I}^{\text{out}}_i(t+1)$ be the internal and external infections at node $i$ be defined as:
\begin{subequations}
\begin{align}
&\bar{I}^{\text{in}}_i(t) \equiv \hiit + \beta_i \hsit\hiit - \gamma_i\hiit + \rho \sum_{j \in \Ni \cap \Vk} \pij \hijt, &&\\    
&\bar{I}^{\text{out}}_i(t) \equiv \rho \sum_{j \in \Ni \cap \Vnk} \pij \hijt. &&
\end{align}
\label{eq:inout}
\end{subequations}

We define the external infection ratio of agent $k$ as: 
\begin{equation}
R_k^{\text{sharing}}(t) \equiv \frac{\sum\limits_{i\in\Vk} \bar{I}^{\text{out}}_i(t)}{\sum\limits_{i\in\Vk}\bar{I}^{\text{in}}_i(t) + \sum\limits_{i\in\Vk} \bar{I}^{\text{out}}_i(t)}.
\end{equation}

At time period $t$, agent $k$ shares $B_k(t) R_k^{\text{sharing}}(t)$ of its budget with other agents. Let $p_{k'k}^I(t)$ be the proportion of infected population traveling from nodes controlled by agent $k'$ to nodes controlled by agent $k$ at time period $t$, i.e.:
\begin{equation}
p_{k'k}^I(t) = \rho \sum_{i \in \Vk}\sum_{j \in \Ni \cap \Vkp} \pij \hijt.
\label{eq:propinf}
\end{equation}

Recall that we assume that agents' budget $B_k(t)$ is function of agents' vaccination capacity $\Gamma_k$. The shared budget of agent $k$ is split among connected agents $k'$ in the global mobility network proportionally to the flow of infected population $p_{k'k}^I(t)$ weighted by the vaccination capacity $\Gamma_{k'}$. Thus, using the proposed budget sharing mechanism, the budget of agent $k$ at time period $t$ is denoted $B^{\text{sharing}}_k(t)$ and is determined as:
\begin{equation}
\label{eq:sharing}
B^{\text{sharing}}_k(t) = B_k(t) (1 - R_k^{\text{sharing}}(t)) + \sum_{k' \in \K \setminus \{k\}} \left(B_{k'}(t) R_{k'}^{\text{sharing}}(t) \frac{\frac{p_{kk'}^I(t)}{\Gamma_{k}}}{\sum\limits_{k'' \in\K} \frac{p_{k''k'}^I(t)}{\Gamma_{k''}}} \right).
\end{equation}

\begin{prop}
The budget sharing mechanism is budget-balanced across all agents, i.e. $\sum_{k \in \K} B_k(t) = \sum_{k \in \K} B^{\text{sharing}}_k(t)$, at any time period $t \in \T$. 
\begin{proof}
Let us rewrite Eq. \eqref{eq:sharing} compactly as $B^{\text{sharing}}_k(t) = B_k(t) (1 - R_k^{\text{sharing}}(t)) + B_k^{\text{received}}(t)$, i.e.:
\begin{equation*}
B_k^{\text{received}}(t) = \sum_{k' \in \K \setminus \{k\}} B_{k'}(t) R_{k'}^{\text{sharing}}(t) \frac{\frac{p_{k k'}^I(t)}{\Gamma_{k}}}{\sum\limits_{k'' \in\K} \frac{p_{k''k'}^I(t)}{\Gamma_{k''}}}    
\end{equation*}

Observe that $B_k(t) (1 - R_k^{\text{sharing}}(t))$ represents the portion of agent' $k$ budget at time $t$ which is not shared with other agents; whereas $B_k^{\text{received}}(t)$ represents the total budget received by agent $k$ from other agents. Since $\sum_{k\in\K} B_k(t) = \sum_{k \in \K} B_k(t) (1 - R_k^{\text{sharing}}(t)) + \sum_{k \in \K} B_k(t) R_k^{\text{sharing}}(t)$, to show that the proposed budget sharing mechanism is budget-balanced, we need only to show that $\sum_{k \in \K} B_k(t) R_k^{\text{sharing}}(t)$ is equal to $\sum_{k\in\K} B_k^{\text{received}}(t)$.

Let $\K = \{1,\ldots,K\}$ be the set of agents. The total budget shared by agents is:
\begin{align}
\sum_{k \in \K} B_k(t) R_k^{\text{sharing}}(t) &= \sum_{k \in \K} B_k(t) R_k^{\text{sharing}}(t) \left(\frac{\frac{p_{1 k}^I(t)}{\Gamma_{1}}}{\sum\limits_{k'' \in\K} \frac{p_{k''k}^I(t)}{\Gamma_{k''}}} + \ldots + \frac{\frac{p_{Kk}^I(t)}{\Gamma_{K}}}{\sum\limits_{k'' \in\K} \frac{p_{k''k}^I(t)}{\Gamma_{k''}}} \right) \nonumber \\
&= B_1(t) R_1^{\text{sharing}}(t)\frac{\frac{p_{1 1}^I(t)}{\Gamma_{1}}}{\sum\limits_{k'' \in\K} \frac{p_{k''1}^I(t)}{\Gamma_{k''}}} + \ldots + B_1(t) R_1^{\text{sharing}}(t)\frac{\frac{p_{K1}^I(t)}{\Gamma_{K}}}{\sum\limits_{k'' \in\K} \frac{p_{k''1}^I(t)}{\Gamma_{k''}}} \nonumber \\
&+ \ldots + B_K(t) R_K^{\text{sharing}}(t)\frac{\frac{p_{1 K}^I(t)}{\Gamma_{1}}}{\sum\limits_{k'' \in\K} \frac{p_{k''K}^I(t)}{\Gamma_{k''}}} + \ldots + B_K(t) R_K^{\text{sharing}}(t)\frac{\frac{p_{KK}^I(t)}{\Gamma_{K}}}{\sum\limits_{k'' \in\K} \frac{p_{k''K}^I(t)}{\Gamma_{k''}}} \nonumber \\ 
&= B_1(t) R_1^{\text{sharing}}(t)\frac{\frac{p_{1 1}^I(t)}{\Gamma_{1}}}{\sum\limits_{k'' \in\K} \frac{p_{k''1}^I(t)}{\Gamma_{k''}}} + \ldots + B_K(t) R_K^{\text{sharing}}(t)\frac{\frac{p_{1 K}^I(t)}{\Gamma_{1}}}{\sum\limits_{k'' \in\K} \frac{p_{k''K}^I(t)}{\Gamma_{k''}}} \nonumber \\
&+ \ldots + B_1(t) R_1^{\text{sharing}}(t)\frac{\frac{p_{K1}^I(t)}{\Gamma_{K}}}{\sum\limits_{k'' \in\K} \frac{p_{k''1}^I(t)}{\Gamma_{k''}}} + \ldots + B_K(t) R_K^{\text{sharing}}(t)\frac{\frac{p_{KK}^I(t)}{\Gamma_{K}}}{\sum\limits_{k'' \in\K} \frac{p_{k''K}^I(t)}{\Gamma_{k''}}} \label{eq:last}
\end{align}

Observe that $p_{kk}^I(t) = 0$ for any agent $k$, hence Eq. \eqref{eq:last} can be rewritten as:
\begin{align*}
\sum_{k \in \K} B_k(t) R_k^{\text{sharing}}(t) &= \sum_{k' \in \K \setminus \{1\}} B_{k'}(t) R_{k'}^{\text{sharing}}(t) \frac{\frac{p_{1 k'}^I(t)}{\Gamma_{k}}}{\sum\limits_{k'' \in\K} \frac{p_{k''k'}^I(t)}{\Gamma_{k''}}} \\
&+ \ldots + \sum_{k' \in \K \setminus \{K\}} B_{k'}(t) R_{k'}^{\text{sharing}}(t) \frac{\frac{p_{K k'}^I(t)}{\Gamma_{k}}}{\sum\limits_{k'' \in\K} \frac{p_{k''k'}^I(t)}{\Gamma_{k''}}} \\
&= \sum_{k \in \K} \sum_{k' \in \K \setminus \{k\}} B_{k'}(t) R_{k'}^{\text{sharing}}(t) \frac{\frac{p_{k k'}^I(t)}{\Gamma_{k}}}{\sum\limits_{k'' \in\K} \frac{p_{k''k'}^I(t)}{\Gamma_{k''}}} = \sum_{k \in \K} B_k^{\text{received}}(t).
\end{align*}
\end{proof}
\end{prop}

To implement the proposed budget-balanced resource sharing mechanism, Algorithm \ref{alg:TS} is extended to include, at each time period $t\in\T$ and for each agent $k\in\K$, the computation of $B^{\text{sharing}}_k(t)$ using Eqs. \eqref{eq:inout}-\eqref{eq:sharing}, and by substituting $B_k(t)$ with $B^{\text{sharing}}_k(t)$ in formulation \eqref{eq:pkt}. The pseudo-code of the TS-based vaccine allocation policy with budget sharing is summarized in Algorithm \ref{alg:TS2}.

\begin{algorithm}[h!]
\label{alg:TS2}
\For{$i \in \V$}{
    $a_i \gets 1$\\
    $b_i \gets 1$\\
}
\For{$k \in \K$}{
    $B^{\text{sharing}}_k(0) \gets B_k(0)$\\
}
\For{$t \in \T$}{
    \For{$i \in \V$}{
        $\htit \gets$ sample prior $\mathsf{Beta}(a_i,b_i)$\\
    }
    \For{$k \in \K$}{
        $\xkt \gets$ solve $\hpkt$ using $B^{\text{sharing}}_k(t)$ instead of $B_k(t)$\\
        \For{$i \in \Vk$}{
            $\uxit \gets 1 - \sum_{t' = \max\{t - m_i + 1,1\}}^{t} x_i(t')$\\
        }
    }
    $[\hsitp,\hiitp,\hritp,\hditp]_{i
    \in \V} \gets$ solve \eqref{eq:SIRDx} using $\xt$ and observe $[\tit]_{i\in\V : \xit >0}$ \\
    \For{$i \in \V$}{
        compute $\bar{I}^{\text{in}}_i(t+1)$ and $\bar{I}^{\text{out}}_i(t+1)$ using Eqs. \eqref{eq:inout}\\
    }
    \For{$k \in \K$}{
        $R_k^{\text{sharing}}(t+1) \gets \frac{\sum_{i\in\Vk} \bar{I}^{\text{out}}_i(t+1)}{\sum_{i\in\Vk}\bar{I}^{\text{in}}_i(t+1) + \sum_{i\in\Vk} \bar{I}^{\text{out}}_i(t+1)}$ \\
        \For{$k' \in \K \setminus \{k\}$}{
            $p_{k'k}^I(t+1) \gets \rho \sum_{i \in \Vk}\sum_{j \in \Ni \cap \Vkp} \pij \bar{I}_j(t+1)$\\
        }
    }
    \For{$k \in \K$}{
        compute $B^{\text{sharing}}_k(t+1)$ using Eq. \eqref{eq:sharing}\\
    }    
    \For{$i \in \V$}{
        \If{$\xit > 0$}{
            $y_i \gets$ Bernoulli trial $(\tit)$\\
            \If{$y_i = 1$}{
                $a_i \gets a_i + 1$ \\
            }
            \Else{
                $b_i \gets b_i + 1$ \\
            }
        }
    }
}
\caption{TS-based policy for online vaccine allocation with budget-balanced resource sharing}
\end{algorithm}

\section{Data}
\label{data}

We conduct a global computational study to test the proposed vaccine allocation policies and explore the potential benefits of the proposed budget sharing mechanism. For this study, we generate a global model of the world wherein agents represent the main countries worldwide. We first describe the population, flight and epidemic data used in this study in Section \ref{popdata} and introduce vaccination related data in Section \ref{capdata}. 

\subsection{Population, flight and epidemic data}
\label{popdata}

World population data was obtained from the Socioeconomic Data and Applications Center \citep{SEDAC} in the form of a uniform 50 km x 50 km raster model. Country boundaries were obtained from GADM \citep{GADM}. Raster population data was converted into individual nodes, forming the set $\V$ of nodes in the network and linked to countries. A total of $|\K|=177$ countries are considered in this study. To implement the radiation model, we assume a distance threshold $D=100$ km for ground mobility. We use a commute fraction  $\alpha_i=0.11$ for all nodes $i \in \V$. A similar approach has been used in previous studies \citep{van2011gleamviz}. Global flight data was processed to extract airports that contribute the most towards the global flight mobility. Airport location data was generated from OurAirports \citep{ourairports}, and the global flight schedule data corresponding to one week of air traffic in October 2020 was obtained from Cirium \citep{Cirium}. Over 392,000 flight entries between the 25$^{th}$ of October to the 31$^{st}$ of October 2020 were analyzed. This was considered a typical sample representing weekly worldwide flights that take place during the pandemic. Inbound and outbound airports in the flight data were used to generate a Voronoi tessellation of the world. Each node of set $\V$ was assigned to an airport based on the obtained Voronoi tessellation. The resulting global network has a total of 53,445 nodes and 12,112,618 mobility links -- including both ground and air travel. 

\begin{figure}
    \centering
    \includegraphics[width=14.5cm]{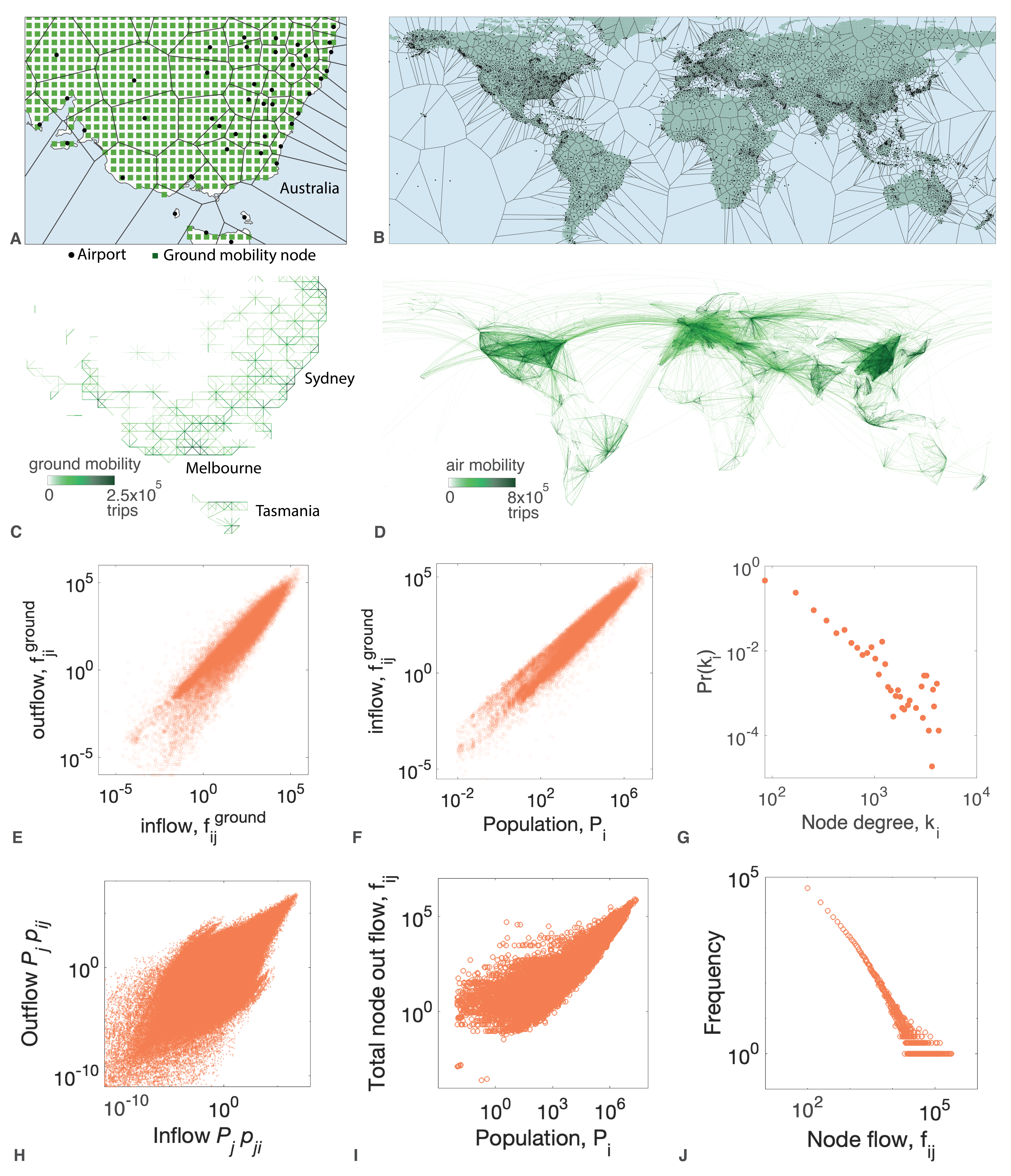}
    \caption{Global human mobility network model. (A) Example representation of airports, cells and Voronoi polygons. (B) The uniform raster model of the world consisting of 50 km by 50 km cells and the Voronoi polygons associated with airports. (C) Ground mobility network consisting of nodes and arcs where arc weights represent ground flow between each pair of nodes estimated with a radiation model. (D) Worldwide air mobility network. (E) The estimated ground mobility flow matrix is asymmetric with high correlation between inflows and outflows for each node. (F) the ground mobility flow is proportional to node population. (G) Near power law node degree distribution of the global human mobility network with combined ground and air flows. (H) Symmetric relationship between node inflow and outflow. (I) Total node outflow $[f_{ij}]_{i \in \V, j \in \Ni}$ is proportional to node population $[P_i]_{i \in \V}$, and (J) distribution of node outflow $[f_{ij}]_{i \in \V, j \in \Ni}$.}
    \label{fig:1}
\end{figure}

The resulting global mobility network and its main features are summarized in Figure 
\ref{fig:1}. Figs. \ref{fig:1}(A) and \ref{fig:1}(B) illustrate the world raster model and the obtained Voronoi tesselation. Fig. \ref{fig:1}(C) illustrate the radiation model used to represent ground mobility and Fig. \ref{fig:1}(D) depicts the flight network used to generate the air mobility component of the global network. Fig. \ref{fig:1} (E) shows the relationship between ground mobility inflow and outflows, while Fig. \ref{fig:1} (F) shows the relationship between ground inflows and node populations. This analysis reveals that the ground mobility matrix is asymmetric and that mobility flows are proportional to nodes population. Fig. \ref{fig:1}(G) illustrates the global mobility network model and reveals a near power law node degree distribution. Fig. \ref{fig:1}(H) shows that inflow and outflow population-weighted mobility rates are nearly symmetric, and Fig. \ref{fig:1}(I) shows that node ouflow are globally proportional to node populations. This shows that the proposed global mobility model meets the conditions required by \citet{brockmann2013hidden} (see Section \ref{epi}). Fig~\ref{fig:1}(J) illustrates the distribution of node outflows in the global human mobility network model. 

Country-based data for COVID-19 transmission and recovery rates and initial susceptible and infected populations are taken from \citet{abbott2020estimating,abbott2020epiforecasts}. Node-based transmission and recovery rates ($\beta_i$ and $\gamma_i$, for all $i \in \V$) and initial susceptible and infected populations are assumed to be uniform for each country. Node-based case fatality rates are set to $\lambda_i = 0.01$ which is equivalent to assume that 1\% of infected individuals die from COVID-19 \citep{rajgor2020many}.

\subsection{Vaccination data}
\label{capdata}

To determine the per-period vaccination capacity of agent $k\in\K$ we assume that countries' vaccination capacity is proportional to their Gross Domestic Product (GDP) per capita. We collected weekly vaccination data from USA, UK, and France as reported by local health authorities during March 2021, and used linear regression to estimate countries' vaccination capacities based on their GDP per capita. The per-period vaccination capacity of country $k\in\K$, denoted $\Gamma_k$, is expressed as a percentage of the population that can be vaccinated per unit of time (week). 

Vaccination budgets $B_k(t) = f(\Gamma_k)$ are determined by either of two modes: without budget sharing and with budget sharing. To implement the former, we use the static function $B_k(t) = \Gamma_k \sum_{i \in \Vk} P_i$ which sets the per-period vaccination budget of agent $k$ as a fraction of the total population controlled by agent $k$. To implement the budget-sharing mechanism, we use Eq. \eqref{eq:sharing} to determine $B^{\text{sharing}}_k(t)$ based on $B_k(t)$ (where $B_k(t)$ is determined as in the no-sharing configuration) and substitute $B_k(t)$ by $B^{\text{sharing}}_k(t)$ in formulation \eqref{eq:hpkt}. Using this model, vaccination budgets $B_k(t)$ are expressed in population units, and the per-unit cost to allocate vaccines to node $i\in\V$ is set to $C_i = P_i$ where $P_i$ is the population of node $i$.

We use countries' vaccination capacities to generate node-based vaccine efficiency rates. This is motivated by the observation that countries with a higher vaccination capacity are also more likely to have successful vaccination campaigns \citep{FORMAN2021553}. We assume that mean vaccine efficiency rates are comprised between 0.5 and 0.9 \citep{ijbs59170} and we scale countries' vaccination capacities in this range using a mapping $g(\cdot)$, i.e. $g(\Gamma_k) \in [0.5,0.9]$, for all $k\in\K$. We denote $\epsilon$ the level of uncertainty in vaccine efficiency rates. For each country $k$ and node $i\in\Vk$, we generate the mean vaccine efficiency rate $\theta_i$ by sampling uniformly and randomly in the range $[g(\Gamma_k) - \epsilon, g(\Gamma_k) + \epsilon]$. The vector of mean vaccine efficiency rates $[\theta_i]_{i\in\V}$ is then used to generate problem instances by assuming that the random variables $[\theta_i(t)]_{i\in\V}$ are generated from uniform distributions with support in $[\theta_i - \epsilon, \theta_i + \epsilon]$.

\subsection{Experimental Framework}
\label{exp}

To test the proposed vaccine allocation policies and budget sharing mechanism, we use data described in Sections \ref{popdata} and \ref{capdata} to generate instances of the proposed global vaccine allocation problem. We consider a two-year time horizon and assume that vaccine allocation decisions are made on a weekly basis by all $|\K|=177$ agents. Thus, a total of $|\T|=104$ time periods are modeled. We consider three levels of uncertainty in vaccine efficiency rates: $\epsilon = 10\%$, $20\%$ and $30\%$. For each level of uncertainty, we generate 100 random instances using the process described in Section \ref{capdata} and we report average performance of the proposed vaccine allocation policies over each group of 100 instances. 

We compare the performance of the proposed TS-based policy implemented using Algorithm \ref{alg:TS} with three alternative policies: PB, MA and GY as described in Section \ref{bench}. In addition, we also report the behavior of the system when no vaccination is performed and hereby refer to this scenario as ``No vaccination''. We assume that all agents (i.e. countries) act independently when making vaccine allocation decisions across their populations. In the base case, hereby also referred to as ``no-sharing'', we further assume that there is no sharing of budget among agents. We then compare the outcome of the vaccination policies when the budget-balanced resource sharing mechanism using Algorithm \ref{alg:TS2} is implemented by all agents. 

All algorithms are implemented in Python on a Windows server with 64 Gb of RAM and a processor Core i9 with a CPU of 3.10 GHz. For research and reproduction purposes, all optimization codes and data used in this study are available at the public repository \url{https://github.com/davidrey123/Vaccine_Allocation}.

\section{Numerical Results}
\label{num}

The numerical results are organized as follows: we first analyze the performance of the TS-based vaccine allocation policy (TS) against benchmarking policies (PB, MA and GY) and under varying levels of uncertainty in vaccine efficiency rates in Section \ref{numbench}. We then focus our attention on the proposed TS-based policy and examine the impact of the available budget for vaccine allocation as well as the budget sharing mechanism onto global infections and deaths in Section \ref{numbudget}.

\subsection{Comparison of vaccine allocation policies}
\label{numbench}

\begin{figure}
    \centering
    \includegraphics[width=\textwidth]{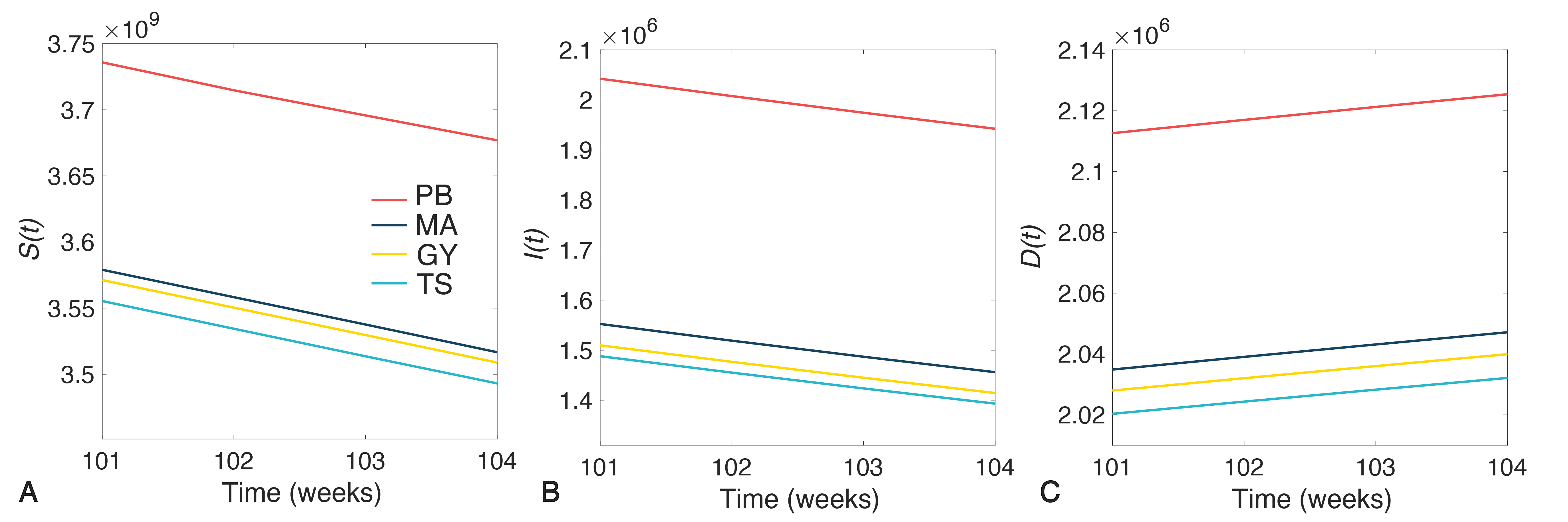}
    \caption{The global number of susceptible, infected, and dead population over the last month of the two-year vaccination planning horizon under different vaccine allocation policies based on population (PB), Moving Average (MA), Greedy (GY), and Thompson Sampling (TS) approaches.}
    \label{fig:2}
\end{figure}

To compare the performance of the proposed TS-based vaccine allocation policy (TS) against benchmarking policies (PB, MA and GY), we implement each policy over 100 instances generated using a level of uncertainty of $\epsilon=20\%$. Fig. \ref{fig:2} reports the average global number of susceptible (S), infected (I) and deaths (D) using each of the four policies over the last month (four time periods) of the two-year vaccination horizon. For all three compartments (S, I and D), the TS-based policy consistently outperforms other policies by achieving a lower number of susceptible, infected and dead population. This suggests that, on average, the TS-based policy learns true mean vaccine efficiency rates faster than the greedy learning policy (GY) which uses the same prior distribution. Both of these Bayesian optimization-based policies outperform the MA policy that only works with historical observations. Compared to the mobility-agnostic PB policy, TS, GY and MA achieves a significantly better performance.

A detailed country-level analysis of the performance of the vaccine allocation policies is summarized in Table \ref{tab:benchmark} where we examine the relative performance of policies MA, GY and TS compared to the PB policy in terms of cumulative and last period gains. The cumulative gain is computed by summing the size of the susceptible population over all 104 time periods of the vaccination horizon whereas the last period gain focuses on the gain achieved in the last week of the horizon. While the cumulative gain is representative of agents' objective function the last period gain better represents the effect of learning throughout the vaccination horizon. Due to space limitations, we focus on the 30 largest countries in terms of number of nodes ($|\Vk|$) which correspond to the optimization problem with the largest number of decision variables. Bold values represent best performance among the three mobility- and disease-aware policies, i.e. MA, GY and TS. Out of 30, the proposed TS-based policy outperforms MA and GY in 22 and 23 cases in terms of cumulative and last period gains, respectively. We observe that for countries in which TS is outperformed by other policies, the performance gap is often marginal while TS is able to substantially improve over MA and GY, such as for Russia (RU), Australia (AU), Mexico (MX) and Indonesia (ID). We find that for Saudi Arabia (SA) the cumulative gains of MA, GY and TS are slightly negative which means that PB outperformed these policies. However, last period gains for SA are positive for MA, GY and TS which suggests that the learning of mean vaccine efficiency rates is gradually improving decision-making. Overall, last period gains tend to be greater than cumulative gains thus reinforcing this hypothesis. A complete table summarizing the performance of the policies over all countries is provided in the public repository linked to this study. This study highlights the role of accounting for global mobility and epidemic dynamics in the design of vaccine allocation policies.\\

\begin{table}
\centering
\begin{tabular}{lrrrrrrrr}
\toprule
& & & \multicolumn{3}{c}{Cumulative Gain (\%)} & \multicolumn{3}{c}{Last Period Gain (\%)} \\
\cmidrule(l){4-6} \cmidrule(l){7-9}
Region & $|\Vk|$ & $\sum_{i \in \Vk} P_i$ & MA & GY & TS & MA & GY & TS \\
\midrule
World & 53445 & 7.28e+09 & 2.21 & 2.26 & \bf 2.40 & 4.36 & 4.57 & \bf 5.00 \\
\midrule
RU & 11604 & 1.42e+08 & 3.98 & 4.04 & \bf 4.20 & 7.34 & 7.58 & \bf 8.31 \\
US & 3877 & 3.34e+08 & 1.35 & 1.37 & \bf 1.57 & 2.64 & \bf 2.74 & 2.31 \\
CN & 3546 & 1.36e+09 & 2.44 & 2.59 & \bf 2.72 & 7.35 & 7.83 & \bf 8.40 \\
BR & 2816 & 1.91e+08 & 5.76 & 5.82 & \bf 5.95 & 12.85 & 13.05 & \bf 13.62 \\
AU & 2755 & 2.28e+07 & 4.67 & 5.33 & \bf 6.02 & 1.34 & 3.94 & \bf 7.87 \\
CA & 2641 & 3.65e+07 & 4.56 & \bf 4.60 & 4.54 & \bf 5.36 & 4.48 & 4.74 \\
KZ & 1323 & 1.91e+07 & \bf 6.36 & 6.29 & 6.14 & \bf 15.21 & 14.98 & 14.56 \\
AR & 1124 & 4.64e+07 & \bf 10.13 & 10.04 & 10.07 & \bf 22.05 & 21.96 & 21.94 \\
IN & 1117 & 1.35e+09 & 0.56 & \bf 0.62 & 0.61 & 1.65 & 1.90 & \bf 1.97 \\
DZ & 855 & 4.09e+07 & 4.10 & 4.02 & \bf 4.17 & 7.90 & 7.95 & \bf 8.29 \\
CD & 749 & 8.71e+07 & 2.62 & 2.62 & \bf 2.64 & 4.46 & 4.47 & \bf 4.62 \\
MN & 741 & 3.28e+06 & 2.19 & 2.33 & \bf 2.44 & 0.12 & 0.56 & \bf 0.98 \\
MX & 703 & 1.27e+08 & 4.93 & 4.95 & \bf 5.34 & 7.23 & 7.14 & \bf 8.36 \\
SA & 683 & 3.30e+07 & -0.88 & -0.85 & \bf -0.79 & 0.68 & 1.41 & \bf 1.99 \\
SD & 630 & 4.53e+07 & \bf 2.16 & 2.12 & 2.14 & 3.76 & 3.77 & \bf 3.85 \\
IR & 621 & 8.09e+07 & 3.72 & \bf 3.75 & \bf 3.75 & 10.64 & 10.64 & \bf 10.81 \\
ID & 599 & 2.30e+08 & 1.32 & 1.52 & \bf 1.79 & 0.02 & 0.61 & \bf 1.55 \\
LY & 592 & 6.58e+06 & 2.56 & 2.64 & \bf 3.02 & 2.07 & 2.34 & \bf 3.55 \\
ZA & 451 & 5.29e+07 & 4.14 & 4.12 & \bf 4.36 & 7.17 & 7.07 & \bf 7.66 \\
TD & 430 & 1.61e+07 & \bf 2.88 & 2.82 & 2.84 & 1.14 & 1.22 & \bf 1.37 \\
PE & 428 & 2.59e+07 & 8.10 & 8.17 & \bf 8.52 & 17.76 & 17.95 & \bf 18.90 \\
ML & 427 & 2.08e+07 & 3.00 & 2.95 & \bf 3.03 & 3.98 & 3.92 & \bf 4.20 \\
AO & 416 & 1.96e+07 & 6.18 & 6.22 & \bf 6.24 & 16.00 & \bf 16.13 & 16.05 \\
NE & 404 & 2.42e+07 & \bf 3.11 & 3.05 & 3.00 & \bf 4.51 & 4.20 & 4.07 \\
CO & 372 & 4.84e+07 & 7.30 & 7.38 & \bf 7.53 & 15.21 & 15.22 & \bf 15.87 \\
ET & 367 & 1.12e+08 & 1.01 & 0.95 & \bf 1.09 & 2.38 & 2.30 & \bf 2.67 \\
BO & 365 & 1.15e+07 & 4.73 & 4.73 & \bf 4.94 & 9.02 & 8.97 & \bf 9.68 \\
MR & 358 & 4.25e+06 & \bf 10.45 & 10.31 & 10.18 & \bf 21.65 & 21.28 & 21.02 \\
TR & 325 & 7.06e+07 & 5.97 & 6.14 & \bf 6.63 & 14.88 & 15.22 & \bf 16.66 \\
PK & 321 & 1.92e+08 & 2.14 & 2.12 & \bf 2.27 & 4.01 & 3.88 & \bf 4.46 \\
\bottomrule
\end{tabular}
\caption{Comparison of vaccine allocation policies for the world and for the 30 countries with the largest number of nodes in the global network. Values reported for policies MA, GY and TS represent the average gain  in percentage over 100 random instances compared to the PB policy. Countries are listed their two-letter ISO code. The Cumulative Gain is computed by summing the size of the susceptible population for each region over all 104 time periods of the vaccination horizon. The Last Period Gain is computed by summing the size of the susceptible population for each region over all 104 time periods of the vaccination horizon. The Last Period Gain corresponds to the gain at the $104^{th}$ time period of the vaccination horizon. Bold values represent best performance.}
\label{tab:benchmark}
\end{table}

\begin{figure}
    \centering
    \includegraphics[width=\textwidth]{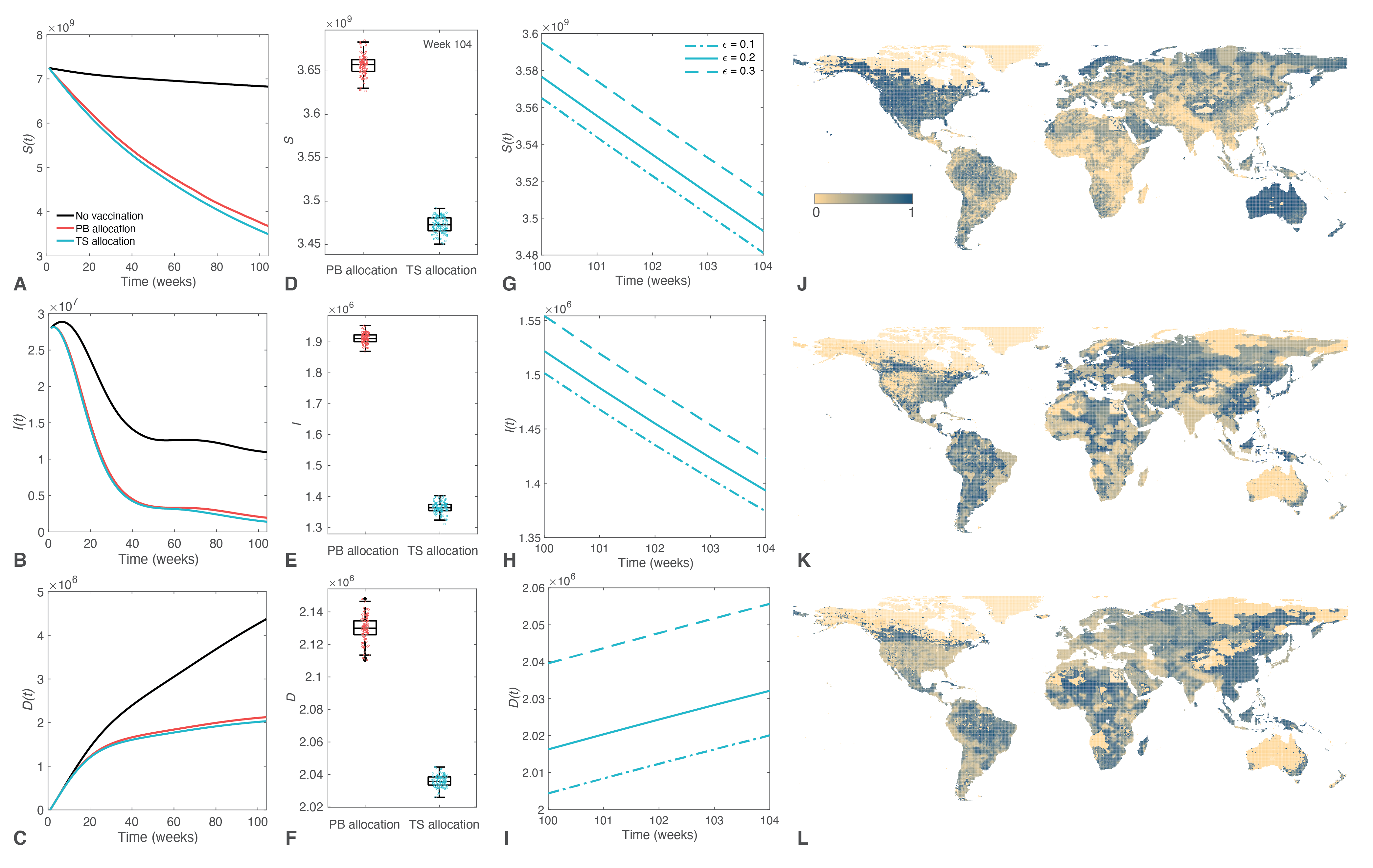}
    \caption{Optimal vaccine allocation to minimize the size of the susceptible population. (A-C) The global number of susceptible, infected, and dead population reduces significantly under a basic population-based vaccine allocation strategy and the proposed reinforcement learning based vaccine allocation strategy over a two-year planning horizon. (D-F) A significantly larger reduction in the size of the susceptible, infected, and dead population is achieved at the end of two-year horizon under the proposed reinforcement learning based vaccine allocation policy. (G-I) Impact of the uncertainty of the vaccine efficiency rate with a focus on the last four weeks of the two-year horizon. (J-L) The global distribution of the percentage reduction in the number susceptible, infected, and dead population across all nodes when a TS-based policy is implemented compared with the base case without vaccination.}
    \label{fig:3}
\end{figure}

In the remaining, we focus on examining the performance of the proposed the TS-based vaccine allocation policy compared to the PB policy that overlooks the impact of mobility and infection transmission dynamics. Fig. \ref{fig:3} provides a comprehensive summary of the global impact of the proposed policies. Figs. \ref{fig:3}(A-C) depicts the average evolution of the susceptible, infected and dead populations, respectively, over the entire two-year vaccination horizon under the TS and PB policies, and in the "No vaccination" scenario. These trends represent average population values across all 100 instances tested using an uncertainty level of $\epsilon= 20\%$. Figs. \ref{fig:3}(D-F) focus on the last week of the two-year horizon and show the distribution of these population values over all 100 instances used in the study. The size of the susceptible and dead population in the world at the end of the two-year horizon with the TS-based vaccine allocation policy is 5\% and 4.4\% smaller, respectively than the size of the susceptible and dead population when a population-based vaccine allocation is used. The impact of the TS-based vaccine allocation policy on the infected population is even larger with 28.3\% reduction in the size compared to that obtained using the PB policy. The impact of the level of uncertainty in the vaccine efficiency rate ($\epsilon$) is examined in Figs. \ref{fig:3}(G-I). Reducing the uncertainty of the vaccine efficiency rate from 20\% to 10\% further reduces the size of the susceptible and dead population by 0.34\% and 0.59\%, while increasing this level of uncertainty to 30\% increases these populations by 0.55\% and 1.16\%, respectively. 

Figs. \ref{fig:3}(J-L) depict the spatial distribution of the resulting reduction in the number of infections and deaths across the world using the TS-based policy with a level of uncertainty of $\epsilon= 20\%$ compared to the no-vaccination case. We find that this reduction is heterogeneous given the heterogeneity in the infection transmission dynamics and human mobility patterns. Interestingly, vast areas in countries such as U.S.A., Canada, Norway, Iceland, Saudi Arabia, and Australia experience significant percentage reductions in the size of their susceptible population while many areas in India, China, Bangladesh, Myanmar, Thailand, Ethiopia, and Nigeria only experience slight percentage reductions. We believe that this is mainly due to the difference in the initial size of susceptible populations, as well as the allocation budget and vaccine administration capacities that are both assumed to be dependent on the GDP per capita of each country. However, in terms of the percentage reduction in death, we observe a significant percentage reduction in vast areas of China, India, Thailand, Vietnam, Indonesia, Brazil, South Africa, and Ethiopia; while many areas in countries including U.S.A., U.K., France, Spain, Iran, and Turkey experience a relatively smaller percentage reduction in death. The observed differences are due to the complex interdependencies of the population density, human mobility, GDP per capita, and infection transmission dynamics.

\subsection{Budget availability and sharing analysis}
\label{numbudget}

\begin{figure}
    \centering
    \includegraphics[width=\textwidth]{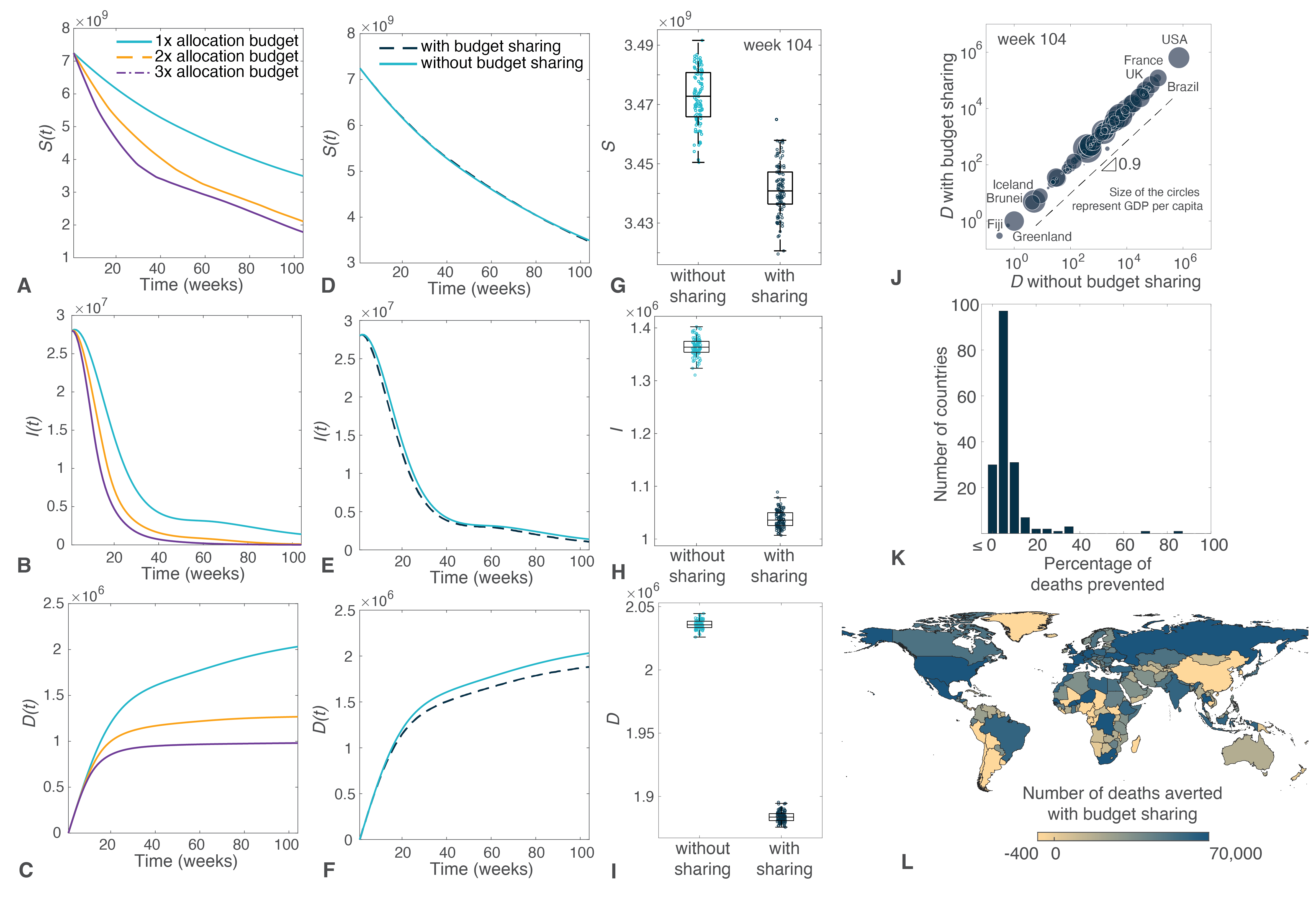}
    \caption{The impact of budget sharing on the world-wide number of infections and deaths. (A-C) The reduction in the global number of susceptible, infected, and dead population when countries’ vaccine allocation budgets are made two and three times higher. (D-F) The impact of budget sharing on the global number of susceptible, infected, and dead population. (G-I) A significant reduction in the global number of infected and dead is achieved via the proposed budget-balanced resource sharing mechanism while the total susceptible population is slightly reduced given a total fixed allocation budget. (J) The number of deaths prevented per individual country under the budget sharing mechanism at the end of the two-year horizon. (K and L) The frequency distribution of the percentage of deaths prevented and the global spatial distribution of the number of deaths prevented across countries under the proposed budget sharing mechanism.}
    \label{fig:4}
\end{figure}

In this section, we study the global behavior of the epidemic under varying vaccination budgets. For this analysis, we focus on the TS-based vaccine allocation policy and set the level of uncertainty in mean vaccine efficiency rates to $20\%$. We implement the budget-balanced resource sharing mechanism proposed in Section \ref{sharing} and compares its performance against a no-sharing mechanism. We also implement the TS-based policy with twofold and threefold inflated vaccine allocation budgets. These numerical results are summarized in Fig. \ref{fig:4}. The impact of the proposed budget sharing mechanism compared to a ``no sharing'' approach is depicted in Figs. \ref{fig:4}(A-F). Under the proposed budget sharing mechanism, while the global size of the susceptible population remains roughly unchanged (see Figs. \ref{fig:4}(A and D)), the size of the infected and dead population is reduced by 24\% and 7.5\% at the end of the two-year vaccination horizon, respectively, compared to an allocation scheme without budget sharing (see Figs. \ref{fig:4}(B, C, E and F)). This averts more than 150,000 deaths and 327,000 new infections over the two-year study horizon. Figs. \ref{fig:4}(J-L) provide further details on the distribution of benefits obtained using the proposed budget sharing mechanism. Notably, we find that the vast majority of countries may substantially reduce their national death toll using this budget sharing mechanism.

We also show that the impact of increasing the global vaccine allocation budget follows a non-linear trend as observed in Figs. \ref{fig:4}(G-I). While a twofold increase of the vaccine allocation budget reduces the size of the susceptible and dead populations by 40\% and 38\%, respectively, further increasing the allocation budget to threefold results in relatively smaller additional gains. Our findings suggest that through cooperation between countries the proposed budget sharing mechanism provides a reduction in the number of deaths equivalent to an increase of 12\% of the global allocation budget.

\section{Conclusion}
\label{con}

This study addressed the problem of allocating vaccines for epidemic control. We consider population subgroups representative of spatial regions connected via a global mobility network and use a compartmental epidemic model to capture disease dynamics. We propose a data-driven optimization approach to solve this vaccine allocation problem in an online fashion. We take the perspective of decision-making agents that aim to minimize the size of their susceptible population and must allocate vaccines under limited supply represented by a budget constraint. We assume that vaccine efficiency rates are unknown and that agents learn these rates from past vaccine allocation decisions. We develop a learning and optimization approach based on Thompson sampling (TS) to learn mean vaccine efficiency rates over time. We propose a budget-balanced resource sharing mechanism to promote cooperation among agents by tracking the source of infections within the global mobility network. 

To explore the behavior of the proposed vaccine allocation policy and mechanisms, we apply the proposed framework to the COVID-19 pandemic. We conduct a global study using a raster model of the world where agents represent the main countries worldwide and have limited vaccine supply. Using real population, flight and epidemic data, we construct a global mobility network that combines both ground and air mobility flows and generate multiple random vaccine allocation problem instances over a two-year vaccination horizon. We then implement the proposed TS-based vaccine allocation policy and benchmark its performance against a population-based (PB) policy, as well as a moving average (MA) policy and a greedy learning (GY) policy. To promote research and result reproduction all optimization codes and data used in this study are made available on a public repository linked to this study (see Section \ref{exp}).

Our numerical results reveal that on average the proposed TS-based policy outperforms the three benchmark policies and leads to reduced susceptible populations as well as lower global number of infections and deaths. Furthermore, our analysis shows that global cooperation in governance and allocation of COVID-19 vaccines could not only reduce worldwide infections and deaths, but also benefit most countries due to the crucial role of human mobility in the spreading of infectious diseases. Notably, countries that have a high mobility exchange can significantly benefit from pooling and sharing their resources. This calls for a more integrated health care management paradigm across policy-makers.

The application of the proposed vaccine allocation policy and the revenue-neutral sharing mechanism to real-world mobility and epidemic data suggests that the proposed methods are of practical use at the global scale. Nevertheless, several modeling assumptions and data sources could be refined to improve the global model, especially in regions where data availability from public health authorities is poor. Future research will explore the use of more detailed epidemic compartmental model available in the literature to improve the accuracy of disease spreading dynamics, e.g. using age-stratified population subgroups. The modeling of the impact of vaccine allocation decisions could be refined by incorporating additional features such as competition effects in the vaccine market \citep{martonosi2021pricing}. The analysis of the impact of coalition among agents could also be investigated to develop further incentive mechanisms to improve global vaccination efforts in the context of a pandemic.

\bibliographystyle{spr-chicago}
\bibliography{scibib}


\end{document}